\documentclass{article}
\usepackage{amssymb}
\usepackage{amsfonts}
\usepackage{amsmath}

\newtheorem{theorem}{Theorem}

\newtheorem{lemma}[theorem]{Lemma}

\newtheorem{remark}[theorem]{Remark}

\newenvironment{proof}[1][Proof]{\noindent\textbf{#1.} }{\
\rule{0.5em}{0.5em}}
\include {mak}
\begin{document}

\title{On some properties on bivariate Fibonacci and Lucas polynomials}
\author{Hac\`{e}ne Belbachir \ and \ Farid Bencherif\ \ \  \and \ \ \ \ \ \
\ \ \ \ \ \ \ \ \ \ \ \ \ \ \ \ \  \\
USTHB, Faculty of Mathematics,\\
Po. Box 32, El Alia,16111, Algiers, Algeria.\\
hbelbachir@usthb.dz or hacenebelbachir@gmail.com\\
fbencherif@usthb.dz or fbencherif@gmail.com}
\maketitle

\begin{abstract}
In this paper we generalize to bivariate polynomials of Fibonacci and Lucas,
properties obtained for Chebyshev polynomials. We prove that the coordinates
of the bivariate polynomials over appropriate basis are families of integers
satisfying remarkable recurrence relations.
\end{abstract}

\section{Introduction}

In \cite{bel10}, the authors established that Chebyshev polynomials of the
first and second kind admit remarkable integer coordinates on specific
basis. It turns out that this property can be extended to Jacobsthal
polynomials \cite{djo3, hor2}, Vieta polynomials \cite{2vie, 2hor12, 2rob1,
2sha1}, and Morgan-Voyce polynomials \cite{2mor, 2and3, 2hor4, 2lee, 2and1,
2swa2, 2che} and Quasi-Morgan-Voyce polynomials \cite{2hor2}, and more
generally to bivariate polynomials associated to recurrence sequences of
order two.

The bivariate polynomials of Fibonacci and Lucas, denoted respectively by $%
(U_{n})=(U_{n}(x,y))$ and $(V_{n})=(V_{n}(x,y))$, are polynomials belonging
to $\mathbb{Z}[x,y]$ and defined by%
\begin{equation*}
\left\{
\begin{array}{l}
U_{0}=0,\text{ }U_{1}=1, \\
U_{n}=xU_{n-1}+yU_{n-2}\text{, }\left( n\geq 2\right)%
\end{array}%
\right. \ \ \ \text{ and}\ \ \ \left\{
\begin{array}{l}
V_{0}=2,\text{ }V_{1}=x, \\
V_{n}=xV_{n-1}+yV_{n-2}\text{, }\left( n\geq 2\right)%
\end{array}%
\right.
\end{equation*}

It is established, see for example \cite{bel1}, that%
\begin{eqnarray}
U_{n+1} &=&\sum_{k=0}^{\left[ n/2\right] }{\binom{n-k}{k}}x^{n-2k}y^{k}, \\
V_{n} &=&\sum_{k=0}^{\left[ n/2\right] }\frac{n}{n-k}{\binom{n-k}{k}}%
x^{n-2k}y^{k}.
\end{eqnarray}

Let $\mathcal{E}_{n}$ be the $%
\mathbb{Q}
$-vectorial space spanned by the free family $\mathcal{C}%
_{n}=(x^{n-2k}y^{k})_{k},$ $0\leq k\leq \left\lfloor n/2\right\rfloor $.
Thus the relations $(1)$ and $(2)$ appear as the decompositions of $U_{n+1}$
and $V_{n}$\ over the canonical basis $\mathcal{C}_{n}$ of $\mathcal{E}_{n}$.

Let us set%
\begin{eqnarray*}
\mathcal{B}_{U} &=&\mathcal{B}_{n,U}=(x^{n-k}U_{n+k+1})_{0\leq k\leq n} \\
\mathcal{B}_{V} &=&\mathcal{B}_{n,V}=(x^{n-k}V_{n+k})_{0\leq k\leq n} \\
\mathcal{B}_{U}^{\ast } &=&\mathcal{B}_{n,U}^{\star
}=(x^{n-k}U_{n+k})_{0\leq k\leq n-1} \\
\mathcal{B}_{V}^{\star } &=&\mathcal{B}_{n,V}^{\star
}=(x^{n-k}V_{n+k-1})_{0\leq k\leq n-1}.
\end{eqnarray*}

The goal of this paper is to prove that for $n\geq 1$, $\mathcal{B}_{U}$\
and $\mathcal{B}_{V}$\ (resp. $\mathcal{B}_{U}^{\ast }$\ and $\mathcal{B}%
_{V}^{\star }$) are basis of $E_{2n}$ (resp. $E_{2n-1}$)\ with respect to
which, the polynomials $U_{2n+1}$ and $V_{2n}$ (resp. $U_{2n}$ and $V_{2n-1}$%
) admit remarkable integer coordinates.

\section{Main results}

\begin{theorem}
We have the following results

\begin{enumerate}
\item $\mathcal{B}_{n,U}$ and $\mathcal{B}_{n,V}$ are basis of $E_{2n},$

\item $\mathcal{B}_{n,U}^{\star }$ and $\mathcal{B}_{n,V}^{\star }$ are
basis of $E_{2n-1}.$
\end{enumerate}
\end{theorem}

As $U_{n+1}$ and $V_{n}$ belong to $E_{n}$, the polynomials $U_{2n+1}$ and $%
V_{2n}$ are elements of $E_{2n}$ with basis $\mathcal{B}_{U}$ or $\mathcal{B}%
_{V}.$ Similarly, $U_{2n}$ and $V_{2n-1}$ belong to $E_{2n-1}$ with basis $%
\mathcal{B}_{U}^{\star }$ or $\mathcal{B}_{V}^{\star }$.

Therefore, there are a priori 8 possible decompositions:%
\begin{equation*}
\begin{array}{lllll}
& \text{over }\mathcal{B}_{U}\text{ }%
\begin{tabular}{|c|}
\hline
1 \\ \hline
\end{tabular}%
\rightarrow trivial, &  &  & \text{over }\mathcal{B}_{U}\text{ }%
\begin{tabular}{|c|}
\hline
3 \\ \hline
\end{tabular}%
\rightarrow simple, \\
U_{2n+1}%
\begin{array}{l}
\nearrow \\
\searrow%
\end{array}
&  &  & V_{2n}%
\begin{array}{l}
\nearrow \\
\searrow%
\end{array}
&  \\
& \text{over }\mathcal{B}_{V}\text{ }%
\begin{tabular}{|c|}
\hline
2 \\ \hline
\end{tabular}%
\rightarrow Th.\ C, &  &  & \text{over }\mathcal{B}_{V}\text{ }%
\begin{tabular}{|c|}
\hline
4 \\ \hline
\end{tabular}%
\rightarrow trivial, \\
&  & \ \ \ \ \ \  &  &  \\
& \text{over }\mathcal{B}_{U}^{\star }\text{ }%
\begin{tabular}{|c|}
\hline
5 \\ \hline
\end{tabular}%
\rightarrow Th.\ A, &  &  & \text{over }\mathcal{B}_{U}^{\star }\text{ }%
\begin{tabular}{|c|}
\hline
7 \\ \hline
\end{tabular}%
\rightarrow Th.\ E, \\
U_{2n}%
\begin{array}{l}
\nearrow \\
\searrow%
\end{array}
&  &  & V_{2n-1}%
\begin{array}{l}
\nearrow \\
\searrow%
\end{array}
&  \\
& \text{over }\mathcal{B}_{V}^{\star }\text{ }%
\begin{tabular}{|c|}
\hline
6 \\ \hline
\end{tabular}%
\rightarrow Th.\ D, &  &  & \text{over }\mathcal{B}_{V}^{\star }\text{ }%
\begin{tabular}{|c|}
\hline
8 \\ \hline
\end{tabular}%
\rightarrow Th.\ B,%
\end{array}%
\end{equation*}%
where the cases 1 and 4 are obvious since $U_{2n+1}\in \mathcal{B}_{U}$ and $%
V_{2n}\in \mathcal{B}_{V}$.

The decomposition of $V_{2n}$ in $\mathcal{B}_{U}$ is simple: we have $%
V_{2n}=2U_{2n+1}-xU_{2n}$.

\ \ \ \ \ \ \ \ \ \ \ \ \ \ \ \ \

Tthe remaining cases are established by the five following results.

\begin{theorem}
(\textbf{A}). Decomposition of $U_{2n+1}$ on basis $\mathcal{B}_{V}.$

For every integer $n\geq 0$, one has%
\begin{equation*}
2U_{2n+1}=\sum_{k=0}^{n}a_{n,k}x^{n-k}V_{n+k},
\end{equation*}%
where%
\begin{equation*}
a_{n,k}=(-1)^{k+1}\binom{n}{k}+2(-1)^{n-k}\sum_{j=0}^{n}(-1)^{j}{\binom{j}{%
n-k}}.
\end{equation*}

Moreover, $\left( c_{n,k}\right) _{n,k\geq 0}$ is a family of integers
satisfying the following recurrence relations%
\begin{equation*}
\left\{
\begin{array}{l}
a_{n,k}=a_{n-1,k}-a_{n-1,k-1}+2\delta _{n,k}\text{, \ }(n\geq 1,\text{ }%
k\geq 1) \\
a_{n,0}=1\text{ }\left( n\geq 0\right) \\
a_{0,k}=\delta _{0,k}\text{ }\left( k\geq 0\right) .%
\end{array}%
\right.
\end{equation*}%
($\delta _{i,j}$ being the Kronecker symbol).
\end{theorem}

The recurrence relations allow obtaining the following table%
\begin{equation*}
\begin{array}{cccccccccc}
\text{\textit{n}}\setminus \text{\textit{k}} & \text{\textit{0}} & \text{%
\textit{1}} & \text{\textit{2}} & \text{\textit{3}} & \text{\textit{4}} &
\text{\textit{5}} & \text{\textit{6}} & \text{\textit{7}} & \text{\textit{8}}
\\
\text{\textit{0}} & 1 &  &  &  &  &  &  &  &  \\
\text{\textit{1}} & 1 & 1 &  &  &  &  &  &  &  \\
\text{\textit{2}} & 1 & 0 & 1 &  &  &  &  &  &  \\
\text{\textit{3}} & 1 & -1 & 1 & 1 &  &  &  &  &  \\
\text{\textit{4}} & 1 & -2 & 2 & 0 & 1 &  &  &  &  \\
\text{\textit{5}} & 1 & -3 & 4 & -2 & 1 & 1 &  &  &  \\
\text{\textit{6}} & 1 & -4 & 7 & -6 & 3 & 0 & 1 &  &  \\
\text{\textit{7}} & 1 & -5 & 11 & -13 & 9 & -3 & 1 & 1 &  \\
\text{\textit{8}} & 1 & -6 & 16 & -24 & 22 & -12 & 4 & 0 & 1%
\end{array}%
\end{equation*}%
from which it follows that%
\begin{eqnarray*}
2U_{1} &=&V_{0}, \\
2U_{3} &=&xV_{1}+V_{2}, \\
2U_{5} &=&x^{2}V_{2}+0V_{3}+V_{4}, \\
2U_{7} &=&x^{3}V_{3}-x^{2}V_{4}+xV_{5}+V_{6}.
\end{eqnarray*}

\begin{theorem}
(\textbf{B}). Decomposition of $U_{2n}$ on basis $\mathcal{B}_{U}^{\star }.$

For every integer $n\geq 1$, one has%
\begin{equation*}
U_{2n}=\sum_{k=0}^{n-1}b_{n,k}x^{n-k}U_{n+k}
\end{equation*}%
where%
\begin{equation*}
b_{n,k}=(-1)^{n-k+1}{\binom{n}{k}}.
\end{equation*}

Moreover, $\left( b_{n,k}\right) _{n,k\geq 0}$ is a family of integers
satisfying recurrence relations%
\begin{equation*}
\left\{
\begin{array}{l}
b_{n,k}=-b_{n-1,k}+b_{n-1,k-1}\text{ }(n\geq 1,\text{ }k\geq 1), \\
b_{n,0}=(-1)^{n+1}\text{ }\left( n\geq 0\right) , \\
b_{0,k}=-\delta _{0,k}\text{ }(k\geq 0).%
\end{array}%
\right.
\end{equation*}
\end{theorem}

The latter recurrence relations allow obtaining the following table%
\begin{equation*}
\begin{array}{ccccccc}
\text{\textit{n}}\setminus \text{\textit{k}} & \text{\textit{0}} & \text{%
\textit{1}} & \text{\textit{2}} & \text{\textit{3}} & \text{\textit{4}} &
\text{\textit{5}} \\
\text{\textit{0}} & -1 &  &  &  &  &  \\
\text{\textit{1}} & \ \ 1 & -1 &  &  &  &  \\
\text{\textit{2}} & -1 & \ \ 2 & -1 &  &  &  \\
\text{\textit{3}} & \ \ 1 & -3 & \ \ 3 & -1 &  &  \\
\text{\textit{4}} & -1 & \ \ 4 & -6 & \ \ 4 & -1 &  \\
\text{\textit{5}} & \ \ 1 & -5 & \ \ 10 & -10 & \ \ 5 & -1%
\end{array}%
\end{equation*}

from which it follows that%
\begin{eqnarray*}
U_{2} &=&\ \ xU_{1} \\
U_{4} &=&-x^{2}U_{2}+2xU_{3} \\
U_{6} &=&x^{3}U_{3}-3x^{2}U_{4}+3xU_{5} \\
U_{8} &=&-x^{4}U_{4}+4x^{3}U_{5}-6x^{2}U_{6}+4xU_{7}
\end{eqnarray*}

\begin{theorem}
(\textbf{C}). Decomposition of $V_{2n-1}$ on basis $\mathcal{B}_{U}^{\star
}. $

For every integer $n\geq 1,$ one has
\begin{equation*}
V_{2n-1}=\sum_{k=0}^{n-1}c_{n,k}x^{n-k}U_{n+k}\text{ with }c_{n,k}=2\left(
-1\right) ^{n-k+1}\dbinom{n}{k}-\delta _{n-1,k}
\end{equation*}

Moreover, $\left( c_{n,k}\right) _{n,k\geq 1}$ is a family of integers
satisfying recurrence relations
\begin{eqnarray*}
c_{n,k} &=&-c_{n-1,k}+c_{n-1,k-1}-\delta _{n,k+2}\text{, \ }(n\geq 2\text{, }%
k\geq 1) \\
c_{n,0} &=&2\left( -1\right) ^{n+1}-\delta _{n,1}\text{ \ }\left( n\geq
1\right) \\
c_{1,k} &=&\left\{
\begin{array}{c}
1\text{ \ \ if \ }k=0 \\
-2\text{\ \ if \ }k=1 \\
0\text{ \ \ if \ }k\geq 2%
\end{array}%
\right. \text{ }
\end{eqnarray*}
\end{theorem}

The latter recurrence relations allow obtaining the following table%
\begin{equation*}
\begin{array}{ccccccc}
\text{\textit{n}}\setminus \text{\textit{k}} & \text{\textit{0}} & \text{%
\textit{1}} & \text{\textit{2}} & \text{\textit{3}} & \text{\textit{4}} &
\text{\textit{5}} \\
\text{\textit{1}} & \ \ 1 & -2 &  &  &  &  \\
\text{\textit{2}} & -2 & \ \ 3 & -2 &  &  &  \\
\text{\textit{3}} & \ \ 2 & -6 & \ \ 5 & -2 &  &  \\
\text{\textit{4}} & -2 & \ \ 8 & -12 & \ \ 7 & -2 &  \\
\text{\textit{5}} & \ \ 2 & -10 & \ \ 20 & -20 & \ \ 9 & -2%
\end{array}%
\end{equation*}

from which we get
\begin{equation*}
\left\{
\begin{array}{l}
V_{1}=\ \ xU_{1} \\
V_{3}=-2x^{2}U_{2}+3xU_{3} \\
V_{5}=\ \ 2x^{3}U_{3}-6x^{2}U_{4}+5xU_{5} \\
V_{7}=-2x^{4}U_{4}+8x^{3}U_{5}-12x^{2}U_{6}+7xU_{7}%
\end{array}%
\right.
\end{equation*}

\begin{theorem}
(\textbf{D}). Decomposition of $V_{2n-1}$ on basis $B_{V}^{\ast }.$

For every integer $n\geq 1,$ one has
\begin{equation*}
2V_{2n-1}=\sum_{k=0}^{n-1}d_{n,k}x^{n-k}V_{n+k-1}\text{ and }d_{n,k}=\left(
-1\right) ^{n-k+1}\frac{n+k}{n}{\binom{n}{k}}
\end{equation*}
\end{theorem}

Moreover,$\left( d_{n,k}\right) _{n\geq 1,\ k\geq 0}$ is a family of
integers satisfying recurrence relations
\begin{eqnarray*}
d_{n,k} &=&-d_{n-1,k}+d_{n-1,k-1}\text{, \ }(n\geq 2\text{, }k\geq 1) \\
d_{n,0} &=&\left( -1\right) ^{n+1}\text{ \ }\left( n\geq 1\right) \\
d_{1,k} &=&\left\{
\begin{array}{c}
1\text{ \ \ if \ }k=0 \\
-2\text{\ \ if \ }k=1 \\
0\text{ \ \ if \ }k\geq 2%
\end{array}%
\right. \text{ }
\end{eqnarray*}

where recurrence relations allow obtaining the following table%
\begin{equation*}
\begin{array}{cccccccc}
\text{\textit{n}}\setminus \text{\textit{k}} & \text{\textit{0}} & \text{%
\textit{1}} & \text{\textit{2}} & \text{\textit{3}} & \text{\textit{4}} &
\text{\textit{5}} & \text{\textit{6}} \\
\text{\textit{1}} & \text{ \ }1 & -2 &  &  &  &  &  \\
\text{\textit{2}} & -1 & \text{ \ }3 & -2 &  &  &  &  \\
\text{\textit{3}} & \text{ \ }1 & -4 & \text{ \ }5 & -2 &  &  &  \\
\text{\textit{4}} & -1 & \text{ \ }5 & -9 & \text{ \ }7 & -2 &  &  \\
\text{\textit{5}} & \text{ \ }1 & -6 & \text{ \ }14 & -16 & \text{ \ }9 & -2
&  \\
\text{\textit{6}} & -1 & \text{ \ }7 & -20 & \text{ \ }30 & -25 & \text{ \ }%
11 & -2%
\end{array}%
\end{equation*}

from which we obtain
\begin{equation*}
\left\{
\begin{array}{l}
2V_{1}=\ \ xV_{0} \\
2V_{3}=-x^{2}V_{1}+3xV_{2} \\
2V_{5}=\ \ x^{3}V_{2}-4x^{2}V_{3}+5xV_{4} \\
2V_{7}=-x^{4}V_{3}+5x^{3}V_{4}-9x^{2}V_{5}+7xV_{6}%
\end{array}%
\right.
\end{equation*}

\begin{theorem}
(\textbf{E}). Decomposition of $U_{2n}$ on basis $\mathcal{B}_{V}^{\star }.$

For every integer $n\geq 1,$ one has%
\begin{equation*}
2U_{2n}=\sum_{k=0}^{n-1}e_{n,k}x^{n-k}V_{n+k-1}\text{ with }e_{n,k}=\frac{1}{%
2}\left( a_{n-1,k}+d_{n,k}\right) +\delta _{n,k}.
\end{equation*}

Moreover, $\left( c_{n,k}\right) _{n,k\geq 0}$ is a family of integers
satisfying recurrence relations%
\begin{equation*}
e_{n,k}=-e_{n-1,k}+e_{n-1,k-1}+a_{n-1,k}\text{ \ }(n\geq 2\text{, }k\geq 1)%
\text{,}
\end{equation*}%
with $e_{n,0}=\left( 1-\left( -1\right) ^{n}\right) /2$ $\ \left( n\geq
1\right) ,\ $and $\ e_{1,k}=\delta _{0,k}$ \ $(k\geq 0).$
\end{theorem}

The latter recurrence relations allow obtaining the following table%
\begin{equation*}
\begin{array}{ccccccc}
\text{\textit{n}}\setminus \text{\textit{k}} & \text{\textit{0}} & \text{%
\textit{1}} & \text{\textit{2}} & \text{\textit{3}} & \text{\textit{4}} &
\text{\textit{5}} \\
\text{\textit{1}} & 1 &  &  &  &  &  \\
\text{\textit{2}} & 0 & \ \ 2 &  &  &  &  \\
\text{\textit{3}} & 1 & -2 & \ \ 3 &  &  &  \\
\text{\textit{4}} & 0 & \ \ 2 & -4 & \ \ 4 &  &  \\
\text{\textit{5}} & 1 & -4 & \ \ 8 & -8 & \ \ 5 &  \\
\text{\textit{6}} & 0 & \ \ 2 & -8 & \ \ 14 & -12 & 6%
\end{array}%
\end{equation*}

from which, we have
\begin{equation*}
\left\{
\begin{array}{l}
2U_{2}=\ \ xV_{0} \\
2U_{4}=0x^{2}V_{0}+2xV_{2} \\
2U_{6}=\ \ x^{3}V_{2}-2x^{2}V_{3}+3xV_{4} \\
2U_{8}=0x^{4}V_{3}+2x^{3}V_{4}-4x^{2}V_{5}+4xV_{6}%
\end{array}%
\right.
\end{equation*}

\section{Proof of Theorems}

Theorem 1 follows from the following lemma.

\begin{lemma}
$\det_{\mathcal{C}_{2n}}($ $\mathcal{B}_{n,U})=\det_{\mathcal{C}_{2n-1}}(%
\mathcal{B}_{n,U}^{\star })=1$ \ and \ $\det_{\mathcal{C}_{2n}}($ $\mathcal{B%
}_{n,V})=\det_{\mathcal{C}_{2n-1}}($ $\mathcal{B}_{n,V}^{\star })=2.$
\end{lemma}

\begin{proof}
Let us prove only the first equality as the proofs of the other ones are
similar.%
\begin{eqnarray*}
\det\nolimits_{\mathcal{C}_{2n}}(\mathcal{B}_{n,U}) &=&\det\nolimits_{%
\mathcal{C}_{2n}}(W_{0},W_{1},...,W_{n})\text{ \ \ \ \ where }%
W_{k}=x^{n-k}U_{n+k+1} \\
&=&\det\nolimits_{\mathcal{C}%
_{2n}}(W_{0},W_{1}-W_{0},...,W_{n-1}-W_{n-2},W_{n}-W_{n-1}),\text{ }
\end{eqnarray*}%
However,%
\begin{equation*}
W_{j}-W_{j-1}=x^{n-j}U_{n+j+1}-x^{n-j+1}U_{n+j}=x^{n-j}(U_{n+j+1}-xU_{n+j})=x^{n-j}yU_{n+j-1}.
\end{equation*}

Thus,%
\begin{equation*}
\det\nolimits_{\mathcal{C}_{2n}}(\mathcal{B}_{n,U})=\det\nolimits_{\mathcal{C%
}_{2n}}(x^{n}U_{n+1},x^{n-1}yU_{n},x^{n-2}yU_{n+1},...,yU_{2n-1})
\end{equation*}

The "component" of $W_{0}=x^{n}U_{n+1}$ over $x^{2n}$ is equal to $1$.

The "component" of $x^{n-j}yU_{n+j-1}$ over $x^{2n}$ is equal to $0,$ so we
have
\begin{eqnarray*}
\det\nolimits_{\mathcal{C}_{2n}}(\mathcal{B}_{n,U}) &=&\det\nolimits_{%
\mathcal{C}_{2n-2}}(x^{n-1}U_{n},x^{n-2}U_{n+1},...,U_{2n-1}) \\
&=&\det\nolimits_{\mathcal{C}_{2n-2}}(x^{n-1-j}U_{n+j})_{0\leq j\leq n-1} \\
&=&\det\nolimits_{\mathcal{C}_{2n-2}}(\mathcal{B}_{n-1,U}) \\
&=&\det\nolimits_{\mathcal{C}_{2n-4}}(\mathcal{B}_{n-2,U})=...=\det%
\nolimits_{\mathcal{C}_{0}}(\mathcal{B}_{0,U})=1
\end{eqnarray*}
\end{proof}

Let $E,$ $A_{m},$ $B_{m},C_{m},$ $D_{m}$ and $E_{m}$ be the operators of $%
\left( \mathbb{Q}\left[ x,y\right] \right) ^{\mathbb{N}}$ defined by%
\begin{equation*}
E\left( \left( W_{n}\right) _{n}\right) =\left( W_{n+1}\right) _{n}
\end{equation*}%
\begin{eqnarray*}
A_{m} &=&\left( x-E\right) ^{m}+2\sum_{k=1}^{m}E^{k}\left( x-E\right) ^{m-k},%
\text{\qquad }\left( m\geq 0\right) , \\
B_{m} &=&-\left( E-x\right) ^{m},\text{\qquad }\left( m\geq 0\right) , \\
C_{m} &=&2E^{m}+2B_{m}-xE^{m-1},\text{\qquad }\left( m\geq 1\right) , \\
D_{m} &=&\left( E-x\right) ^{m-1}\left( x-2E\right) ,\text{\qquad }\left(
m\geq 1\right) , \\
E_{m} &=&\frac{1}{2}\left( xA_{m-1}+D_{m}+2C_{m}\right) ,\text{\qquad }%
\left( m\geq 0\right) ,
\end{eqnarray*}%
where $E$ is the forward shift operator given by $EW_{n}=W_{n+1}$.

Then, we have%
\begin{eqnarray*}
A_{m} &=&\sum_{k=0}^{m}a_{m,k}x^{m-k}E^{k}\qquad \text{with}\qquad
a_{m,k}=\left( -1\right) ^{k+1}{\binom{m}{k}}+2\left( -1\right)
^{m-k}\sum_{j=0}^{m}\left( -1\right) ^{j}{\binom{j}{m-k}} \\
B_{m} &=&\sum_{k=0}^{m}b_{m,k}x^{m-k}E^{k}\qquad \text{with}\qquad
b_{m,k}=\left( -1\right) ^{m-k+1}{\binom{m}{k}} \\
C_{m} &=&\sum_{k=0}^{m-1}c_{m,k}x^{m-k}E^{k}\qquad \text{with}\qquad
c_{m,k}=2\left( -1\right) ^{m-k+1}{\binom{m}{k}}-\delta _{m-1,k} \\
D_{m} &=&\sum_{k=0}^{m}d_{m,k}x^{m-k}E^{k}\qquad \text{with}\qquad d_{m,k}=%
\frac{\left( -1\right) ^{m-k+1}\left( m+k\right) }{m}{\binom{m}{k}} \\
E_{m} &=&\sum_{k=0}^{m-1}e_{m,k}x^{m-k}E^{k}\qquad \text{with}\qquad e_{m,k}=%
\frac{1}{2}\left( a_{m-1,k}+d_{m,k}\right) +\delta _{m,k}
\end{eqnarray*}

With these notations, relations stated by Theorems $A,$ $B,$ $C,$ $D$ and $E$
may be expressed by means of the following relations%
\begin{eqnarray*}
\mathbf{a}.\text{ }\forall n &\in &\mathbb{N}\ \ \ \ \ A_{n}V_{n}=2U_{2n+1}
\\
\mathbf{b}.\text{ }\forall n &\in &\mathbb{N}\ \ \ \ \ \ B_{n}U_{n}=0 \\
\mathbf{c}.\text{ }\forall n &\in &\mathbb{N}^{\ast }\ \ \ \
C_{n}U_{n}=V_{2n-1} \\
\mathbf{d}.\text{ }\forall n &\in &\mathbb{N}^{\ast }\ \ \ \ D_{n}V_{n-1}=0
\\
\mathbf{e}.\text{ }\forall n &\in &\mathbb{N}^{\ast }\ \ \ \
E_{n}V_{n-1}=2U_{n}
\end{eqnarray*}

which are to be proven. For this, the following lemma will be useful for us.

\begin{lemma}
For every integers $n$ and $m,$ we have

\begin{enumerate}
\item $\left( x-E\right) ^{n}U_{m}=\left( -y\right) ^{n}U_{m-n}$ \ and $\
\left( x-E\right) ^{n}V_{m}=\left( -y\right) ^{n}V_{m-n}$ \ $\left( m\geq
n\geq 0\right) $

\item $V_{n}=2U_{n+1}-xU_{n}$ \ $\left( n\geq 0\right) $

\item $V_{n}=U_{n+1}+yU_{n-1}$ \ $\left( n\geq 1\right) $

\item $\sum_{k=1}^{n}\left( -y\right) ^{n-k}V_{2k}=U_{2n+1}-\left( -y\right)
^{n}$ \ \ $\left( n\geq 0\right) $
\end{enumerate}
\end{lemma}

\begin{proof}
\ \ \ \ \ \ \ \ \ \ \ \ \ \ \ \ \ \ \ \ \ \ \ \ \ \ \ \ \

\begin{enumerate}
\item We proceed by induction on $n$, observing that for $n=1$, we have $%
\left( x-E\right) ^{n}U_{m}=xU_{m}-U_{m+1}=-yU_{m-1}$ and $\left( x-E\right)
^{n}V_{m}=-yU_{m-1}$, for $m\geq 1.$

\item For every integer $n\in \mathbb{N}$, let us put $%
S_{n}:=2U_{n+1}-xU_{n} $. We observe that $S_{0}=2$, $S_{1}=x$ and $%
S_{n}=xS_{n-1}+yS_{n-2}$ for $n\geq 2.$ Thus, for every $n\in \mathbb{N}$, $%
V_{n}=S_{n}=2U_{n+1}-xU_{n}.$

\item For every integer $n\geq 1$, we have from the latter relation $%
V_{n}=U_{n+1}+\left( U_{n+1}-xU_{n}\right) =U_{n+1}+yU_{n-1}.$

\item For every integer $n\in \mathbb{N}$, put $T_{n}:=U_{2n+1}-%
\sum_{k=1}^{n}\left( -y\right) ^{n-k}V_{2k}.$ The relation to be proven is
equivalent to $T_{n}=\left( -y\right) ^{n}$ \ $\left( n\geq 0\right) .$
Then, we remark that from relation 1. of this lemma, we have for every
integer $n\geq 1$%
\begin{equation*}
T_{n}+yT_{n-1}=U_{2n+1}+yU_{2n-1}-V_{2n}=0
\end{equation*}%
$\left( T_{n}\right) _{n\geq 0}$ is then a geometric sequence with
multiplier $\left( -y\right) $ and of first term $T_{0}=1$. It follows that
for every integer $n\in \mathbb{N}$, $T_{n}=\left( -y\right) ^{n}.$
\end{enumerate}
\end{proof}

\begin{proof}[Proof of relations a., b., c., d. and e]
\ \ \ \ \ \ \ \ \ \ \ \ \ \ \ \ \ \ \ \ \ \ \ \ \ \ \ \ \ \ \ \ \ \ \ \ \

\textbf{a}. For every integer $n\in \mathbb{N}$, we have

$\qquad \qquad \qquad \qquad \qquad \qquad
\begin{array}{ll}
A_{n}V_{n} & =(\left( x-E\right) ^{n}+2\sum_{k=1}^{n}E^{k}\left( x-E\right)
^{n-k})V_{n} \\
& =\left( -y\right) ^{n}V_{0}+2\sum_{k=1}^{n}\left( -y\right) ^{n-k}V_{2k}%
\text{ \ \ \ \ \ \ (from 1 of lemma)} \\
& =2U_{2n+1}%
\end{array}%
$

\textbf{b}. For every integer $n\in \mathbb{N}$, we have

$\qquad \qquad \qquad \qquad \qquad \qquad
\begin{array}{ll}
B_{n}U_{n} & =-\left( E-x\right) ^{n}U_{n} \\
& =-\left( -y\right) ^{n}U_{0}\text{ \ \ \ \ (from 1 of lemma)} \\
& =0%
\end{array}%
$

\textbf{c}. For every integer $n\in \mathbb{N}^{\ast }$, we have

$\qquad \qquad \qquad \qquad \qquad \qquad
\begin{array}{ll}
C_{n}U_{n} & =\left( 2E^{n}+2B_{n}-xE^{n-1}\right) U_{n} \\
& =2U_{2n}+2B_{n}U_{n}-xU_{2n-1},%
\end{array}%
$

\ \ however $B_{n}U_{n}=0$ (from a.), thus

$\qquad \qquad \qquad \qquad \qquad \qquad
\begin{array}{ll}
C_{n}U_{n} & =2U_{2n}-xU_{2n-1} \\
& =V_{2n-1}%
\end{array}%
$

\textbf{d}. For every integer $n\in \mathbb{N}^{\ast }$, we have

$\qquad \qquad \qquad \qquad \qquad \qquad
\begin{array}{ll}
D_{n}V_{n-1} & =\left( x-2E\right) \left( E-x\right) V_{n-1} \\
& =\left( x-2E\right) V_{0} \\
& =xV_{0}-2V_{1} \\
& =0%
\end{array}%
$

\textbf{e}. For every integer $n\in \mathbb{N}^{\ast }$, we have

$\qquad \qquad \qquad \qquad \qquad \qquad
\begin{array}{ll}
E_{n}V_{n-1} & =\left( \frac{1}{2}xA_{n-1}+\frac{1}{2}D_{n}+E^{n}\right)
V_{n-1}\medskip \\
& =\frac{1}{2}xA_{n-1}V_{n-1}+\frac{1}{2}D_{n}V_{n-1}+V_{2n-1}\medskip%
\end{array}%
$

\qquad But $A_{n-1}V_{n-1}=2U_{2n-1}$ (from c.) and $D_{n}V_{n-1}=0$ (from
b.), It follows that

$\qquad \qquad \qquad \qquad \qquad \qquad
\begin{array}{ll}
E_{n}V_{n-1} & =xU_{2n-1}+V_{2n-1} \\
& =2U_{2n}\text{ \ \ \ \ (From 2 of the lemma).}%
\end{array}%
$
\end{proof}

\begin{remark}
Theorems A, B, C, D and E generalize results obtained for the Chebyshev
polynomials \cite{bel10}, Indeed,
\begin{eqnarray*}
\frac{1}{2}V_{n}(2x,1) &=&T_{n}(x)\text{ is the Chebyshev polynomials of the
first kind, } \\
U_{n+1}(2x,1) &=&U_{n}(x)\text{ is the Chebyshev polynomials of the second
kind,}
\end{eqnarray*}%
with%
\begin{equation*}
\left\{
\begin{array}{l}
T_{n}(x)=2xT_{n-1}-T_{n-2} \\
T_{0}=1,T_{1}=x%
\end{array}%
\right. \text{ \ and \ }\left\{
\begin{array}{l}
U_{n}(x)=2xU_{n-1}-U_{n-2} \\
U_{0}(x)=1,U_{1}=2x%
\end{array}%
\right.
\end{equation*}
\end{remark}

\end{document}